\documentclass[11pt]{amsart} 
\usepackage{latexsym, amssymb, bbm}
\usepackage{amsthm}
\usepackage[all]{xy}

\newtheorem{lem}{Lemma}
\newtheorem{thm}[lem]{Theorem}
\newtheorem{prop}[lem]{Proposition}
\newtheorem{cor}[lem]{Corollary}

\theoremstyle{definition}
\newtheorem{defn}[lem]{Definition}
\newtheorem{example}[lem]{Example}

\newtheorem{remark}[lem]{Remark}

\DeclareMathOperator{\supp}{supp}

\DeclareMathOperator{\vo}{Vol}

\DeclareMathOperator{\sgrad}{sgrad}

\DeclareMathOperator{\Gr}{Gr}
\DeclareMathOperator{\can}{can}

\DeclareMathOperator{\Cal}{Cal}

\newcommand{\wt}[1]{\widetilde{#1}}
\newcommand{\abs}[1]{\left| #1 \right|}
\newcommand{\norm}[1]{\left\| #1 \right\|}
\newcommand{\ind}{\boldsymbol{\mathbbm{1}}}

\def\co{\colon\thinspace}
\def\pc{\!:\!}

\def\C{\mathbb{C}}
\def\N{\mathbb{N}}

\def\R{\mathbb{R}}
\def\Z{\mathbb{Z}}

\def\Si{\Sigma}

\def\a{\alpha}

\def\d{\partial}
\def\Del{\Delta}
\def\e{\epsilon}
\def\i{\iota}

\def\l{\lambda}

\def\r{\rho}
\def\s{\sigma}
\def\t{\tau}
\def\th{\theta}
\def\T{\Theta}

\def\vp{\varphi}

\def\w{\omega}
\def\z{\zeta}
\def\vp{\varphi}
\def\p{\phi}

\def\C{\mathbb{C}}
\def\CP{\mathbb{CP}}

\def\D{\mathbb{D}}

\def\L{\mathcal{L}}
\def\La{\Lambda}

\def\O{\Omega}
\def\H{\mathcal{H}}
\def\P{\mathcal{P}}
\def\PH{\P Ham}
\def\CH{\widetilde{Ham}}

\newcommand{\param}[1]{\left(#1\right)}

\begin{document}
\title[Symplectic reduction of quasi-morphisms and quasi-states]{Symplectic reduction of\\ quasi-morphisms and quasi-states}
\author{Matthew Strom Borman}
\address{Department of Mathematics, University of Chicago, Chicago, Illinois 60637}
\email{borman@math.uchicago.edu}
\thanks{This work is partially supported by the NSF-grant DMS 1006610.}

\begin{abstract}
We prove that quasi-morphisms and quasi-states on a closed integral symplectic manifold 
descend under symplectic reduction to symplectic hyperplane sections.
Along the way we show that quasi-morphisms that arise from spectral invariants
are the Calabi homomorphism when restricted to Hamiltonians supported on stably displaceable
sets.
\end{abstract}  

\maketitle


\section{Introduction and results}

\subsection{An overview}
In the series of papers 
\cite{EntPol03, EntPol06, EntPol08, EntPol09RS}, Entov and Polterovich introduced a way to
construct quasi-morphisms on the universal cover of the Hamiltonian group
and symplectic quasi-states, when $(M, \w)$ is a closed symplectic manifold.
In order to determine if their construction will work for a specific $(M^{2n}, \w)$, one must
compute the quantum homology ring $QH_{2n}(M,\w)$ and determine if it contains a field summand,
which in practice can be a nontrivial task.
In this paper we introduce a reduction method
that allows one to build quasi-morphisms and symplectic quasi-states from known examples, without
requiring further quantum homology calculations.  This process is presented formally in 
Theorems~\ref{thm: descent qm} and \ref{thm: descent}, which are the main results of this paper.

As an example, 
suppose that $(M, \w; \Si)$ forms a subcritical polarization where $\Si \subset M$ is a complex
hypersurface.  If $\z\co C^\infty(M) \to \R$ is a symplectic quasi-state that
vanishes on functions with displaceable support, then our construction produces a symplectic
quasi-state $\bar{\z}\co C^\infty(\Si) \to \R$.  The main tool of the construction is
Biran's decomposition theorem \cite{Bir01}, which roughly states that 
the unit disk normal bundle $E_\Si$, for the hypersurface $\Si \subset M$, is symplectomorphic to
$M \setminus \Del$, where in the subcritical case $\Del$ is a displaceable set (see Section~\ref{Biran} for a precise statement).  This allows us to define a map
$\T\co C^\infty(\Si) \to C^\infty(M)$, which lifts functions $H$ on $\Si$ to 
$E_\Si$,
multiplies them by a function that vanished at the boundary of the fibers, and then views the
result as a function on $M$ using Biran's symplectomorphism.  The reduced symplectic
quasi-state $\bar{\z}\co C^\infty(\Si) \to \R$ is the result of pulling $\z$ back by the map $\T$.
The construction for quasi-morphisms proceeds in a similar fashion.

This reduction method produces quasi-morphisms and symplectic quasi-states on the Hirzebruch surfaces, which were shown to have quasi-morphisms and quasi-states by Ostrover \cite{Ost06}, but required quantum homology calculations.  These and other examples are discussed in Sections~\ref{subcritical} and \ref{critical}.

\subsection{The Hamiltonian group}
Let $(M^{2n}, \w)$ be a closed symplectic manifold, any smooth function
$H \in C^\infty(M)$ determines a Hamiltonian vector field $\sgrad H$ on $M$ by
$\w(\sgrad H, \cdot) = -dH$.  In this manner any time-dependent Hamiltonian $F\co M \times [0,1] \to \R$ gives a Hamiltonian isotopy $\phi_F = \{f_t\}_{t \in [0,1]}$, by solving the differential equation
$\d_t f_t = (\sgrad F_t)_{f_t}$ with initial condition $f_0 = \ind$.  The collection of all such time one maps $f_1$ is the \emph{Hamiltonian group} $Ham(M,\w)$.  Seen as a Lie group, its Lie algebra is $\H(M)$,
all $H \in C^\infty(M)$ that are \emph{normalized} to have mean zero $\int_M H \w^n = 0$,
and the bracket is the Poisson bracket
\[
	\{H, K\} = \w(\sgrad K, \sgrad H) = dH(\sgrad K).
\]

Any smooth path in $Ham(M)$ is a Hamiltonian isotopy, so the space of smooth paths based
at the identity $\PH(M)$, can be identified with $\P\H(M)$, the space of functions
$F\co M \times [0,1] \to \R$ such that $F_t \in \H(M)$ at all times.  
The group structure of time-wise product on $\PH(M)$ carries over to $\P\H(M)$ as
$\phi_F\phi_G = \phi_{F\#G}$ and $\phi_F^{-1} = \phi_{\bar{F}}$
where
\[
	(F\#G)(x,t) = F(x,t) + G(f_t^{-1}(x),t) \quad \mbox{and} \quad \bar{F}(x,t) = -F(f_t(x),t).
\]
The universal cover $\CH(M)$ is the path space $\PH(M)$ where elements are considered up to
homotopy with fixed endpoints.  See \cite{McDSal98, Pol01} for further preliminaries and note we 
use the sign conventions of \cite{Pol01}.

\subsection{Quasi-morphisms on $\CH(M)$}

A \emph{homogeneous quasi-morphism} on a group $G$ is a function $\mu\co G \to \R$ so that 
$n\mu(g) = \mu(g^n)$ for all $n \in \Z$ and $g \in G$, and for some $C \geq 0$:
\begin{equation}\label{quasi}
	\abs{\mu(g_1 g_2) - \mu(g_1) - \mu(g_2)} \leq C \quad \mbox{for all $g_1, g_2 \in G$}.
\end{equation}
The smallest such $C$ is the \emph{defect} and is denoted $D(\mu)$.  Homogeneous quasi-morphisms
are conjugation invariant and they are genuine homomorphisms when restricted to abelian subgroups.  On perfect groups, like $Ham(M)$ and $\CH(M)$, they are the best one can hope for in terms of a map to 
$\R$.  If $\mu$ only statisfies \eqref{quasi}, then $\mu$ is a \emph{quasi-morphism} and can be homogenized by forming $\check{\mu}(g) = \displaystyle{\lim_{n \to \infty}\tfrac{\mu(g^n)}{n}}$, which will be the unique homogeneous quasi-morphism within a bounded distance from $\mu$.  See 
\cite{Cal09, Kot04} for more information about quasi-morphisms.

In \cite{EntPol03}, Entov and Polterovich introduced a way to build homogeneous quasi-morphisms
$\mu\co\CH(M) \to \R$ using spectral invariants from Hamiltonian Floer theory.  We will briefly recall
the construction and refer the reader to \cite{EntPol03, EntPol08, EntPol09RS} for more details.
Let $QH_{*}(M,\w)$ be the quantum homology ring for $(M, \w)$, which is 
$H_*(M; \C) \otimes_\C \La_\w$
where $\La_\w$ is a Novikov ring associated to $(M,\w)$, and its ring structure is given
by a deformation of the normal intersection product using Gromov-Witten invariants.  To each element
$a \in QH_*(M)$ there is an associated spectral invariant, defined in terms of Hamiltonian Floer theory,
which is a functional $c(a, \cdot)\co C^\infty(M \times [0,1]) \to \R$.

With $c(a, \cdot)$, one can define the functional $\mu(a, \cdot)\co C^\infty(M \times [0,1]) \to \R$:
\begin{equation}\label{qm def}
	\mu(a, F) = \int_0^1\int_M F(x,t)\, \w^n dt -\vo(M) \lim_{k\to \infty} \frac{c(a,F^{\#k})}{k},
\end{equation}
which descends to a function
\begin{equation}
	\mu(a, \cdot)\co \CH(M) \to \R.
\end{equation}
If $a \in QH_{2n}(M)$ is an idempotent and $a * QH_{2n}(M)$ is a field \cite[Remark 1.23]{EntPol09RS},
then $\mu(a, \cdot)$ is a homogeneous quasi-morphism.  Recently Usher \cite[Theorem 1.4]{Ush10} and Fukaya, Oh, Ohta, and Ono \cite[Theorem 1.2]{FukOhOht11a} have independently extended this construction using idempotents in the big quantum homology ring, and the associated spectral invariants.  Any quasi-morphism build using spectral invariants we will call a \emph{spectral quasi-morphism}. 

Spectral quasi-morphisms have the following two additional
properties, where $\CH_U(M)$ are all Hamiltonian paths generated by Hamiltonians with support in $U$
(not necessarily normalized):
\begin{enumerate}
	\item \emph{Calabi Property:} If $U \subset M$ is open and displaceable, then
	$\mu$ restricted to $\CH_U(M)$ is the Calabi homomorphism $\Cal_U$ given by
	\[\Cal_U(\phi_F) = \int_0^1\int_U F_t\, \w^n dt.
	\]
	\item \emph{Stability Property:} There is some $B > 0$, so that for all $F, G \in \P\H(M)$,
	\[
		\int_0^1 \min_M(F_t - G_t)\,dt \leq \frac{\mu(\phi_G) - \mu(\phi_F)}{B} \leq
		\int_0^1 \max_M(F_t-G_t)\,dt.
	\]
	In the case of spectral quasi-morphisms $B$ can be taken to be $\vo(M)$.  We will
	say that a quasi-morphism $\mu$ on $\CH(M)$ is \emph{stable} if it has the stability
	property.
\end{enumerate}
A subset $X \subset M$ is \emph{displaceable} if there is a $\vp \in Ham(M,\w)$ so that
$\vp(X) \cap \overline{X} = \emptyset$ and $X$ is \emph{stably displaceable} if
$X \times S^1 \subset M \times T^*S^1$ is displaceable.
Stable displaceability does not imply displaceability as shown in
\cite[Example 1.30]{EntPol09RS}.  For spectral quasi-morphisms, the Calabi property can be strengthen
to include stably displaceable sets, which will be proved in Section~\ref{calabi}.

\begin{thm}\label{thm: strong calabi}
	Let $U \subset (M, \w)$ be a stably displaceable open subset, then
	spectral quasi-morphisms $\mu(a, \cdot)\co \CH(M) \to \R$ restrict to the
	Calabi homomorphism $\Cal_U$ on $\CH_U(M)$.
\end{thm}

Recently Shelukhin \cite[Corollary 1]{She11} proved that there is a quasi-morphism on $\CH(M)$ for any closed symplectic manifold, using a construction that is not Floer-theoretic. 
Theorem~\ref{thm: descent qm} and \ref{thm: descent} below cannot be directly applied to Shelukhin's quasi-morphisms since they do not have the Calabi property and they do not induce quasi-states.

\subsection{Symplectic quasi-states and quasi-measures}
\emph{Symplectic quasi-stats} are functionals 
$\z\co C^\infty(M) \to \R$ such that for $H,K \in C^\infty(M)$, $a \in \R$:
	\begin{enumerate}
		\item \emph{Quasi-linearity}: If $\{H, K\} = 0$, then $\z(H + aK) = \z(H) + a\,\z(K)$.
		\item \emph{Monotonicity}: If $H \leq K$, then $\z(H) \leq \z(K)$.
		\item \emph{Normalization}: $\z(1) = 1$.
	\end{enumerate}
Note that these properties imply that $\z$ is Lipschitz with respect to the $C^0$-norm,
$\abs{\z(H) - \z(K)} \leq \norm{H-K}$.

This notion was introduced by Entov and Polterovich \cite{EntPol06} and is the symplectic version of Aarnes' notion of a \emph{topological quasi-state} \cite{Aar91}.  These are functionals $\z\co C(X) \to \R$, for $X$ a compact space, that are monotonic, normalized, and quasi-linear in the sense that 
if $F \in C(X)$ and $g,h\in C(\R)$, then $\z(g\circ F + h\circ F) = \z(g\circ F) + \z(h\circ F)$.
Aarnes proved a Riesz representation theorem for quasi-states showing they correspond to topological
quasi-measures, which are additive set functions $\t\co \mathcal{A} \to [0,1]$ where $\mathcal{A}$ are the subsets of $X$ that are either open or closed.  The quasi-state $\z$ determines the quasi-measure
$\t_\z$, where for a closed set $C \subset X$
\[
\t_\z(C) = \inf\{\z(F) \mid F: X \to [0,1],\, F \geq \ind_C\}.
\]
Being a symplectic quasi-state is a stronger condition than being a topological quasi-state. 

A stable homogeneous quasi-morphism $\mu \co \CH(M) \to \R$ induces a symplectic quasi-state
$\z_\mu \co C^\infty(M) \to \R$ via
\begin{equation}\label{qs from qm}
	\z_\mu(H) = \frac{\int_M H \w^n - \mu(\phi_{H_n})}{\vo(M)},
\end{equation}
where $H_n = H - \tfrac{\int_M H \w^n}{\vo(M)}$.  
Quasi-linearity of $\z_{\mu}$ follows from the fact that $\p_H$ and $\p_K$ commute if $\{H, K\} = 0$,
and quasi-morphisms are homomorphisms on abelian subgroups.  The monotonicity of
$\z_{\mu}$ follows from the stability property for $\mu$.

Symplectic quasi-states formed as in \eqref{qs from qm} are $Ham(M)$
invariant, since $\mu(\phi_H) = \mu(\psi^{-1} \phi_H \psi)$ by conjugation invariance of quasi-morphisms,
and $\psi^{-1} \phi_H \psi = \phi_{H\circ \psi}$.   If $\mu$ has the Calabi property, then $\z_\mu$ vanishes on functions with displaceable support and hence the associated quasi-measure $\t_\mu$ will vanish on displaceable sets.  Since the symplectic quasi-state $\z_{\mu}$ is $Ham(M)$ invariant, the associated
quasi-measure $\t_\mu$ will also be $Ham(M)$ invariant, 
which allows it to detect nondisplaceability.  For if $\vp\in Ham(M)$ displaced $X \subset M$,
then $\tau_\mu(X) \leq 1/2$ since 
\[2\, \tau_\mu(X) = \tau_\mu(X) + \tau_\mu(\vp(X)) = \tau_\mu(X \cup \vp(X)) \leq 1.\]

Consider now a spectral quasi-morphism $\mu(a, \cdot): \CH(M) \to \R$.  The associated
symplectic quasi-state $\z(a, \cdot): C^\infty(M) \to \R$ is given by
\[
	\z(a, H) = \lim_{k \to \infty} \frac{c(a, kH)}{k},
\]
and we can form the associated symplectic quasi-measure $\t(a, \cdot)$.
In \cite[Definition 1.4]{EntPol09RS}, the notion of a closed subset $X \subset M$ being \emph{(super)heavy} (with respect to an idempotent $a$) was introduced.  When $a$ gives a spectral quasi-morphism, these two notions agree and are equivalent to $\t(a, X) = 1$ \cite[Remark 1.23]{EntPol09RS}.  In this language they prove \cite[Theorem 1.4]{EntPol09RS} that if $X$ is stably displaceable, then $\t(a,X) \not= 1$. This result can be strengthen to the following theorem, which is proven in Section~\ref{proofs2}.

\begin{thm}\label{thm: stably}
	For any subset $V \subset M$, if $V$ is stably displaceable, then $\tau(a, V) = 0$.
\end{thm}

For a symplectic quasi-state $\z\co C^\infty(M) \to \R$, consider the inequality
\begin{equation}\label{master}
	\Pi_\z(H,K) := \abs{\z(H + K) - \z(H) - \z(K)} \leq C \sqrt{\norm{\{H,K\}}}.
\end{equation}
In \cite[Theorem 1.4]{EntPolZap07} it is proved that \eqref{master} holds for
symplectic quasi-states formed from stable homogeneous quasi-morphisms as in \eqref{qs from qm}, where $C$ can be taken to be $\sqrt{2D(\mu)}/\vo(M)$.
If \eqref{master} holds for some $C$, let $C(\z)$ be the smallest number for which it is satisfied
for all $H, K \in C^\infty(M)$.  Due a result of Cardin and Viterbo \cite[Theorem 1.2]{CarVit08}, which was later generalized by Entov and Polterovich \cite[Theorem 1.1]{EntPol10} and
Buhovsky \cite[Theorem 1.1.2]{Buh10}, it makes sense to say that two continuous functions $H$ and $K$ Poisson commute if they have smooth $C^0$-approximations $H_n$ and $K_n$ such that 
$\norm{\{H_n, K_n\}} \to 0$.  Therefore \eqref{master} says that symplectic
quasi-states are linear on continuous functions that Poisson commute.  The inequality \eqref{master} is the main tool used to lower bound the Poisson bracket invariants introduced by Buhovsky, Entov, and Polterovich \cite{BuhEntPol11}.

\subsection{Polarizations and Biran's decomposition theorem}\label{Biran}
\begin{defn}
	A \emph{polarization} $\P=(M^{2n},\O, J;\Si)$ is a closed K\"ahler manifold with 
	$[\O] \in H^2(M,\Z)$ and a closed connected complex
	hypersurface $\Si \subset M$ such that $[\Si] \in H_{2n-2}(M, \Z)$ is Poincar\'e dual
	to $k[\O] \in H^2(M, \Z)$.  By scaling $\O$ by $k$, we will assume that $k = 1$
	and hence $\vo(M, \O) = \vo(\Si, \O_\Si)$.
\end{defn}

Given a polarization $\P= (M, \O, J; \Si)$, there is a compact isotropic celluar subspace
$\Del_{\P} \subset M$, which is disjoint from $\Si$ and is called the \emph{skeleton} associated to 
$\P$.  An example of a polarization is a linear $\CP^{n-1} \subset \CP^n$, where
$\CP^n$ is given the Fubini-Study form $\O$ with $\int_{\CP^1} \O = 1$.  If $\CP^{n-1}$ is given
by $\{z_0 = 0\}$, then the skeleton $\Del = \{[1\pc 0 \pc \cdots \pc 0]\}$ is a single point.

One can build a symplectic disk bundle $\pi\co (E_\Si, \w_{\can}) \to \Si$ using
the normal bundle from $\Si \subset M$, such that $(\Si, \w_\Si)$ symplectically embeds into
$(E_\Si, \w_{\can})$ as the zero section.  
Biran's decomposition theorem for polarizations \cite[Theorem 2.6.A]{Bir01} says there is
	a canonical symplectomorphism $F_\P$ such that the following diagram commutes:
	\begin{equation}\label{decom}
	\xymatrix{
	(E_\Si, \w_{\can}) \ar[r]^{F_\P} & (M \setminus \Del_\P, \O)\\
	(\Si, \O_\Si) \ar[u]^{\mbox{$0$-section}} \ar@{=}[r] & (\Si, \O_\Si) \ar[u]_{\mbox{inclusion}}
	}
	\end{equation}
The disk bundle $E_\Si$ comes with a hermitian metric on the fibers such
that each fiber is the open disk of radius one.  The radius on the fibers
gives a global coordinate function $r$ on $E_\Si$ and the map
$F_\P$ is such that any open neighborhood of $\Del_\P$ contains the complement of
$F_\P(\{e \in E_\Si \mid r(e) \leq 1-\e\})$ for some small $\e$.  For more details see 
Section~\ref{setup} and \cite{Bir01, Bir06, BirCie01}.

\subsection{Symplectic reduction of quasi-morphisms and symplectic quasi-states}
	We will now describe our procedure for reducing quasi-morphisms and symplectic quasi-states
	from the total space $M$ of a polarization to the hypersurface $\Si$, which is the main idea 
	in this paper.
		
		Let $\th\co [0,1] \to \R$ be a nonnegative function such that $\th(0) = 1$, 
		vanishes in a neighborhood of $1$, and 
		defines a smooth function $\th(r)$ on $E_\Si$, where $r$ is the radial coordinate
		on $E_\Si$.
		To $\th$ we associate a linear, order preserving map
		\begin{equation}\label{Theta}
			\T\co C^\infty(\Si) \to C^\infty(M)
			\quad \mbox{given by} \quad 
			\T(H) = (F_\P^{-1})^*(\th \cdot \pi^*H),
		\end{equation}
		where $\pi\co  E_\Si \to \Si$ is the projection and $F_\P$ is the symplectomorphism
		in \eqref{decom}.
		When notationally convenient we will consider
		$\th$ as the the function $\T(1)$ on $M$, and observe that the condition that $\th$ 
		vanishes near $1$ ensures that $\T(H)\co  M \to \R$ is smooth.
		If $H \in C^\infty(\Si)$ has zero mean, then so does $\th\cdot\pi^*H \in C^\infty(E_\Si)$
		(see \eqref{fiber}).
		Therefore we can view $\T$ as a linear map
		\[
			\T\co  \P\H(\Si) \to \P\H(M).
		\]
		A nice class of maps are given by
		$\T_\e$, where
		$\th_\e(r)\co  [0,1] \to \R$ interpolates between
		$1-r^2$ when $r \leq 1-\e$ and zero when $r \geq 1-\tfrac{\e}{2}$.
	
	\begin{thm}\label{thm: descent qm}
		Let $\mu\co  \CH(M) \to \R$ be a homogeneous quasi-morphism
		and 
		suppose there is an open neighborhood $U$ of $\Del_\P$ so that $\mu$ restricts
		to the Calabi homomorphism on $\CH(M)_U$.  
		Then $\mu$ descends to a homogeneous quasi-morphism $\hat{\mu}$, where for
		any $\e$ sufficiently small:
		\begin{equation}\label{qmorph descent}
			\hat\mu\co  \CH(\Si) \to \R \quad \mbox{is defined by} \quad 
			\hat\mu(\vp) = \mu(\phi_{\T_\e(F)}),
		\end{equation}
		where $F \in \P\H(\Si)$ generates
		$\vp \in \CH(\Si)$. 
		The defect of $\hat\mu$ has upper bound $D(\hat\mu) \leq 2D(\mu)$. 
		If $\mu$ is stable, then so is $\hat\mu$, and if $\mu$ also has the Calabi property
		then the normalized version of $\hat\mu$
		\[
			\bar\mu\co  \CH(\Si) \to \R \quad\mbox{defined by}
			\quad \bar\mu(\vp) = \z(\th_\e)^{-1} \mu(\phi_{\T_\e(F)})
		\]
		will also have the Calabi property, where $\z$ is the symplectic quasi-state
		associated to $\mu$.
	\end{thm}

	The condition on $\e$ is that the complement of
	$F_\P(\{e \in E_\Si \mid r(e) \leq 1-\e\})$ is a subset of the neighborhood
	$U$ of $\Del_\P$.  Since $\mu$ is the Calabi homomorphism on $\CH(M)_U$,
	it follows from \eqref{qs from qm} that $\z$ vanishes on functions supported in $U$.
	Therefore $M \setminus U$ has full quasi-measure $\tau_{\z}(M \setminus U) = 1$,
	and hence $\z(\th_{\e}) > 0$ since $\th_{\e}$ is strictly positive on $M \setminus U$.	
	The map $\T$ can also be used to pullback symplectic quasi-states $\z$
	on $M$ to $\Si$.
	
	\begin{thm}\label{thm: descent}
		Suppose that $\z\co  C^\infty(M) \to \R$ is a symplectic quasi-state and let
		$\th$ as in \eqref{Theta} be such that $\z(\th) > 0$.  Then
		the functional $\bar\z_\th\co  C^\infty(\Si) \to \R$, defined as
		\begin{equation}\label{qs reduction}
		\bar\z_\th(F) = \frac{\z(\T(F))}{\z(\th)},
		\end{equation}
		is a symplectic quasi-state.  The properties of vanishing on functions
		with displaceable support and being $Ham$ invariant descend to
		$\bar\z_\th$ should $\z$ have them.  If $\z$ satisfies \eqref{master}, 
		then $\bar\z_\th$ does as well with
		$C(\bar\z_\th) \leq \frac{C(\z)}{\z(\th)} \norm{\tfrac{\th(r)}{\sqrt{1-r^2}}}$.
	\end{thm}
	
	Such $\th$ exist if and only if the skeleton $\Del_{\P}$ does not have
	full quasi-measure with respect to $\t_{\z}$, and in these cases they will exist with abundance.
	Theorems~\ref{thm: descent qm} and \ref{thm: descent} are proved in
	Section~\ref{proofs}, and examples of when they apply 
	are given in Sections~\ref{subcritical} and \ref{critical}.
	
	\begin{remark}	
	If Theorem~\ref{thm: descent qm} applies to
	a stable $\mu:  \CH(M) \to \R$, then we have two ways to build symplectic
	quasi-states on $\Si$.  We can reduce $\mu$ to $\bar\mu$ as in Theorem~\ref{thm: descent qm},
	and then form the associated quasi-state $\z_{\bar\mu}$.  Alternatively we can first form the
	associated quasi-state $\z = \z_\mu$ on $M$, and then reduce it to $\bar\z_{\th_\e}$ as in
	Theorem~\ref{thm: descent}, where $\e$ must be small enough for 
	Theorem~\ref{thm: descent qm} to apply.  These two symplectic quasi-states $\z_{\bar\mu}$ and
	$\bar\z_{\th_\e}$ are equal, which can be verified by checking on normalized Hamiltonians 
	$F \in C^{\infty}(\Si)$
	\[
		\z_{\bar\mu}(F) = -\frac{\bar\mu(\phi_{F})}{\vo(\Si)} 
		= -\frac{\mu(\phi_{\T_{\e}(F)})}{\z(\th_{\e}) \vo(\Si)}
		= \frac{\z(\T_{\e}(F))}{\z(\th_{\e})} = \bar\z_{\th_{\e}}(F).
	\]
	\end{remark}
	
	\begin{remark}
	Theorems \ref{thm: descent qm} and
	\ref{thm: descent} only rely on the Biran decomposition \eqref{decom}, for which the 
	condition that
	$(M, \O, J)$ be K\"ahler is not essential.
	Let $(M^{2n}, \O)$ be a closed symplectic manifold with $[\O]\in H^{2}(M; \Z)$ 
	and let $\Si^{2n-2}$ be a closed symplectic submanifold Poincar\'e dual to $k[\O]$
	for $k \in \N$.  In this setting $(M, \O; \Si)$, Opshtein proved a Biran decomposition result 
	\cite[Theorem 1]{Ops11} where the the skeleton $\Delta$ in \eqref{decom} has zero volume.
	Under the additional assumption that $M \setminus \Si$ is a Weinstein manifold, $\Si$ is called
	a \emph{symplectic hyperplane section} and for these
	Biran and Khanevsky \cite[Section 2.3]{BirKha11} explain that Biran's decomposition holds with
	the skeleton $\Delta$ being an isotropic cellular subspace.	 
	Closed symplectic
	submanifolds $\Si^{2n-2} \subset (M^{2n}, \O)$ that are Poincar\'e dual to $k[\O] \in H^{2}(M; \Z)$ 
	exist for sufficiently large $k \in \N$ by Donaldson \cite{Don96}.  
	For $k$ sufficiently large, these Donaldson 
	symplectic hypersurfaces also satisfy the condition that $M \setminus \Si$ is Weinstein by 
	Giroux \cite[Proposition 11]{Gir02}.
\end{remark}
	
	\begin{remark}
	Suppose that in Theorem~\ref{thm: descent qm} 
	we replaced Biran's symplectomorphism $F_{\P}$ in \eqref{decom} 
	with a symplectic embedding of a disk bundle over 
	$\Si$ into a tubular neighborhood of $\Si$ in $M$.  
	For a small enough tubular neighborhood of $\Si$, 
	reducing any spectral quasi-morphism $\mu(a, \cdot)$ this way
	will result in the trivial quasi-morphism by Corollary~\ref{cor: coiso}.
	Therefore a priori Theorems~\ref{thm: descent qm} and \ref{thm: descent}
	cannot be made into local constructions.
	
	In \cite{Bor11a}, the constructions in this paper are adjusted to the case of symplectic
	reduction for Hamiltonian $\mathbb{T}^{k}$-actions on a level set of the moment
	map $\Phi: (M, \w) \to \R^{k}$ with full 
	quasi-measure.  In contrast to this paper, the constructions in \cite{Bor11a} are truly local in nature.
	\end{remark}
	
	\begin{remark}
	This method of reducing quasi-morphisms and symplectic quasi-states 
	also works in the setting of Liouville domains $(M, \w, L)$, when
	the Reeb vector field for $\a = (\i_L\w)|_{\d M}$ induces a free $S^1$-action on
	$\d M$.  Then the analogues of Theorems~\ref{thm: descent qm} and \ref{thm: descent}
	hold and allow one to reduce quasi-morphism and symplectic quasi-states on $M$ to
	the reduction $\Si = \d M/S^1$.
	\end{remark}	

\subsection*{Acknowledgments}
I would like to thank my advisor Leonid Polterovich for his wonderful help, guidance, and encouragement related to this work.  I am grateful to Paul Biran for his generous help with examples of subcritical polarizations and explaining his work with Yochay Jerby to me.  I would also like to thank Egor Shelukhin for useful discussions.


\section{Remarks and examples}

\subsection{A criterion for a spectral quasi-morphism to restrict to the Calabi homomophism}\label{calabi}

From the definitions, whenever a stable homogeneous quasi-morphism on $\CH(M)$ restricts
to the Calabi homomorphism on $\CH_U(M)$ for an open set $U$, the associated quasi-measure vanishes on $U$.  It turns out that the converse is true as well for spectral quasi-morphisms.

\begin{prop}\label{prop: zero measure}
	Let $\mu(a, \cdot)\co  \CH(M) \to \R$ be a spectral quasi-morphism and let
	$\t(a,\cdot)$ be its associated
	symplectic quasi-measure. If an open subset $U \subset M$ has zero measure $\tau(a,U) =0$,
	then $\mu$ restricts to the Calabi homomorphism on $\CH(M)_U$.
\end{prop}
\begin{proof}
	Let $\vp \in \CH_U(M)$ be generated by a Hamiltonian $F\co  M \times [0,1] \to \R$ such that
	$\supp(F_t) \subset U$ for all $t$.
	The fact that $\tau(a,U) = 0$ implies that $\z(a,H) = 0$ for all $H \in C^\infty(M)$ 
	that are zero outside of $U$.
	There exist functions $A, B\co  M \to \R$ that are zero outside $U$ so that 
	as time dependent functions $kA \leq F^{\#k} \leq kB$ and hence 
	by the monotonicity of spectral invariants
	\[
		0 = \z(A) = \lim_{k \to \infty} \frac{c(a, kA)}{k} \leq 
		\lim_{k \to \infty} \frac{c(a, F^{\#k})}{k} \leq \lim_{k \to \infty} \frac{c(a, kB)}{k}
		= \z(B) = 0.
	\]
	Therefore the limit term vanishes in \eqref{qm def} and hence $\mu(a, \vp) = \Cal_U(\vp)$.
\end{proof}

We can now prove Theorem~\ref{thm: strong calabi}.

\begin{proof}[Proof of Theorem~\ref{thm: strong calabi}]
	By Theorem~\ref{thm: stably}, a stably displaceable set $U$ has zero measure
	and hence by Proposition~\ref{prop: zero measure}, $\mu(a, \cdot)$ restricts to the
	Calabi homomorphism on $\CH_U(M)$.
\end{proof}

A corollary to Theorem~\ref{thm: strong calabi} is the following.

\begin{cor}\label{cor: coiso}
	Let $\Sigma \subset (M, \w)$ be a closed nowhere
	coisotropic submanifold, then $\Sigma$ is stably displaceable and
	hence its has a neighborhood on which $\mu(a, \cdot)$ restricts
	to the Calabi homomorphism.  A special case being
	when $\Sigma \subset M$ is a closed symplectic submanifold.
\end{cor}
\begin{proof}
	Observe that $\Sigma \times S^1 \subset M \times T^*S^1$ is nowhere coisotropic
	and has at least one nonvanishing section of its normal bundle.  Therefore by
	G\"urel's Displacement Principle \cite[Theorem 1.1]{Gur08}, 
	we know that $\Sigma \times S^1$ is displaceable in $M \times T^*S^1$.
\end{proof}

\subsection{Reduction for subcritical polarizations.}\label{subcritical}

The skeleton $\Delta_\P$ of a polarization $\P = (M^{2n}, \O, J; \Si)$ is an
isotropic cellular subspace and can be replaced with an isotropic
CW-complex $\Del$ of the same dimension such that the decomposition \eqref{decom} still holds
 \cite[Theorem 2.6.C]{Bir01}.  A polarization is \emph{subcritical} if $\dim \Del_\P < n$
 and is \emph{critical} if $\dim \Del_\P = n$, which are the only options since $\Del_\P$ is isotropic.
In the subcritical case, the Isotopy Theorem \cite[Theorem 6.1.1]{BirCie01} states that
$\Del$ is displaceable.  The following is now a corollary of Theorems~\ref{thm: descent qm}
and \ref{thm: descent}.

\begin{cor}\label{cor: subcritical}
	Let $\P = (M^{2n},\O, J;\Si)$ be a subcritical polarization.
	Calabi quasi-morphisms $\mu\co  \CH(M) \to \R$ descend to
	$\bar\mu\co  \CH(\Si) \to \R$ as in Theorem~\ref{thm: descent qm}.
	Symplectic quasi-states $\z$ on $(M, \w)$ with the vanishing property
	descend to symplectic quasi-states $\bar{\z}$ on $(\Si, \O_{\Si})$
	as in Theorem~\ref{thm: descent}.
\end{cor}

\subsubsection{The Hirzebruch surfaces are hypersurfaces in subcritical polarizations.}

The Hirzebruch surfaces are given by
\[
	\Si_k = \{(z, w) \in \CP^2 \times \CP^1 \mid z_0w_0^k + z_1 w_1^k = 0\}
\]
for $k \in \N$, when $k$ is even $\Si_k$ is diffeomorphic to $\CP^1 \times \CP^1$,
and when $k$ is odd $\Si_k$ is diffeomorphic to $\CP^2$ blown-up at $[0 \pc 0 \pc 1]$.
Let $\w_n$ be the Fubini-Study form on $\CP^n$ such that $\int_{\CP^1}\w_n = 1$.
Observe that the Poincar\'e dual of $[\Si_k] \in H_4(\CP^2 \times \CP^1;\Z)$ is
$[\O_k] := [\w_2] + k[\w_1]$.  Therefore $\P_k = (\CP^2 \times \CP^1, \O_k, J; \Si_k)$ is a polarization of degree $1$.  

The holomorphic line bundle $\L_k \to \CP^2 \times \CP^1$ determined by $\Si_k$, has an unique (up to scaling) holomorphic section 
$s_k$ with zero set $\Si_k$.  Picking a hermitian metric $\norm{\cdot}$ on $\L_k$ gives a function 
$\r_{k}\co  \CP^2 \times \CP^1 \setminus \Si_k \to \R$ where
\[
	\r_{k}(z,w) = -\norm{s_k(z,w)}^2 
	= -\frac{\abs{z_0w_0^k + z_1 w_1^k}^2}{\norm{z}^2\norm{w}^{2k}}.
\]
The skeleton $\Del_k = \Del_{\P_k}$ is defined to be the union of the stable manifolds for
$\nabla \r_k$.  We have that 
\[
	\Del_k = \{[\bar{w}_0^k \pc \bar{w}_1^k \pc 0] \times [w_0 \pc w_1] \mid w \in \CP^1\}
\]
which is $2$-dimensional, and therefore the polarization $\P_k$ is subcritical.  
See \cite[Section 3]{Bir01} for similar examples.  To summarize:

\begin{prop}\label{main example}
	The polarizations $\P_k= (\CP^2 \times \CP^1, \O_k = \w_2 + k\, \w_1, J; \Si_k)$
	given by the Hizebruch surfaces $\Si_k$ are subcritical and the symplectomorphism
	types of the hypersurfaces are as follows: 
	\begin{enumerate}
	\item $(\Si_{2l}, \O_{2l}|_{\Si_{2l}})$ is symplectomorphic to 
	$(\CP^1 \times \CP^1, \w_1 \oplus 3l\, \w_1)$,
	
	\item $(\Si_{2l+1}, \O_{2l+1}|_{\Si_{2l+1}})$ is symplectomorphic to 
	$(\CP^2\, \#\, \overline{\CP^2}, \w)$, where $\w$ is the symplectic form
	with $\w(L) = 3l+2$ for $L$ the general line and $\w(E) = 3l+1$ for
	$E$ the exceptional divisor.
	\end{enumerate}
	These two families of rational ruled surfaces have quasi-morphisms
	and symplectic quasi-states that descend from $(\CP^2 \times \CP^1, \O_k)$, 
	as in Corollary~\ref{cor: subcritical}.
\end{prop}
\begin{proof}
	The $\Si_k$ are ruled surfaces, since projection onto the $\CP^1$ factor
	$\pi_k\co  \Si_k \to \CP^1$ gives a $\CP^1$-fiber bundle.
	Hence by Lalonde and McDuff \cite{LalMcD96}, the 
	symplectomorphism type of $(\Si_k, \O_{k}|_{\Si_k})$ is determined by the 
	cohomology class $[\O_{k}|_{\Si_k}] \in H^2(\Si_k, \Z)$.
A basis for $H_2(\Si_k, \Z)$ is given by
\[
	F_k= [\pi_k^{-1}([1\pc b]) ]
	= [\{z_0 + z_1b^k = 0\} \times [1\pc b]] \quad \mbox{and} \quad D_k = [[0\pc 0 \pc 1] \times \CP^1],
\]
where $\O_k(F_k) = 1$ and $\O_k(D_k) = k$, and the intersections are 
$F_k \cdot F_k = 0$, $F_k \cdot D_k = 1$, and $D_k \cdot D_k = -k$.

When $k = 2l$ is even, $(\Si_{2l}, \O_{2l}|_{\Si_{2l}})$ and $(\CP^1 \times \CP^1, \w)$ are
symplectomorphic and the homology classes of the spheres are identified with
\[
	[\CP^1 \times pt] = F_{2l} \quad \mbox{and} \quad [pt \times \CP^1] = lF_{2l} + D_{2l}.
\]
Therefore $\w([\CP^1 \times pt]) = 1$ and $\w([pt \times \CP^1]) = 3l$, which proves (1).

When $k= 2l +1$ is odd, $(\Si_{2l+1}, \O_{2l}|_{\Si_{2l}})$ and $(\CP^2\, \#\, \overline{\CP^2}, \w)$
are symplectomorphic and the homology classes of $L$ and $E$ are identified with
\[
L = (l+1)F_{2l+1} + D_{2l+1} \quad \mbox{and}\quad E = l F_{2l+1} + D_{2l+1}.
\]
Therefore $\w(L) = 3l+2$ and $\w(E) = 3l+1$, which proves (2). 
\end{proof}

\subsubsection{Subcritical polarizations from algebraic geometry.}

Biran and Jerby \cite{BirJer11} have proved that if $M \subset \CP^m$ is a smooth algebraic manifold
whose dual variety $M^* \subset \CP^{m*}$ has codimension at least $2$,
then a hyperplane section $\Si$ in $M$ gives a subcritical polarization $(M, \Si)$.
Examples of such $M$ are the complex Grassmannians $\Gr(2, 2n+1)$ and
 $\CP^n$-bundles over smooth projective varieties of dimension less than $n$
 with linear fibers.  Refer to \cite[Section 9.2.C]{Tev03} for other examples.

\subsection{Reduction for critical polarizations.}\label{critical}

If $\P = (M^{2n}, \O, J; \Si, \Del_\P)$ is a critical polarization, then $\Del_\P$ has dimension 
equal to $n$ and is possibly nondisplaceable.  Therefore it is not enough that a quasi-morphism
have the Calabi property for it to descend, however in light of Proposition~\ref{prop: zero measure} we can say the following:

\begin{cor}\label{cor: critical}
	Let $\P = (M^{2n},\w, J;\Si, \Del_\P)$ be a critical polarization.
	Spectral quasi-morphisms $\mu\co  \CH(M) \to \R$ descend to
	$\bar\mu\co  \CH(\Si) \to \R$ as in Theorem~\ref{thm: descent qm}, provided
	$\Del_\P$ has a neighborhood $U$ with zero measure $\t_\mu(U) = 0$.
\end{cor}

\begin{example}\label{EP example}
Let $\Si$ be the diagonal in $(\CP^1 \times \CP^1, \w_1 \oplus \w_1)$.
This forms a critical polarization and the anti-diagonal $\Del$ is the corresponding skeleton
\cite[Section 3.2.1]{Bir01}.  The quantum homology 
$QH_4(\CP^1 \times \CP^1)$ splits into a direct sum of fields, with idempotents
$a_{\pm}$ and corresponding spectral quasi-morphisms $\mu_{\pm}$, quasi-states
$\z_{\pm}$, and quasi-measures $\tau_{\pm}$.  

In \cite[Theorem 1.1]{EliPol10} Eliashberg and Polterovich compute that $\Delta$ has a neighborhood
$U$ with zero measure $\tau_+(U) = 0$, in particular $\tau_+(K) = 1$ where
$K$ is an exotic torus disjoint from $\Del$.
Therefore by Corollary~\ref{cor: critical}, $\mu_+$ and $\z_+$ descend to
$(\Si, \w_1 \oplus \w_1) = (\CP^1, 2\,\w_1)$.
The reduced quasi-state $\bar\z_+$ is the median quasi-state, for it is the unique
symplectic quasi-state that vanishes on displaceable sets \cite[Section 8]{EntPol06},
but the uniqueness of $\mu([S^2], \cdot)$ as a Calabi quasi-morphism on
$Ham(S^2)$ is unknown.
\end{example}

\section{Proofs}

\subsection{The general set-up}\label{setup}

Let $(P, \xi)$ be a closed contact manifold with contact form $\a$ and 
assume the Reeb vector field $R$ generates a free $S^1$-action on $P$.   
The quotient $\Si = P/S^1$, has a symplectic form $\s = -\pi_* d\a$ induced by the projection 
$\pi\co  P \to \Si$.  From this principal $S^1$-bundle, we can build 
the \emph{standard symplectic disk bundle $\pi\co (E, \w) \to (\Si, \s)$ modeled on $P$}.  The
total space $E = P \times_{S^1} \D^2$ has a symplectic
form given by $\w = \pi^*\s + d(r^2\a) = -d((1\!-\!r^2)\a)$, where $r$ is the radial coordinate on the open unit disk 
$\D^2$. Note that $\sgrad r^2 = R$, each fiber has area one (by taking $S^1 = \R/\Z$), and 
$(\Si, \s)$ embeds symplecticly into $E$ as the zero section.
Given a polarization $\P = (M, \O, J; \Si)$, the unit normal bundle 
$P \to \Si$ is a prequantization space for $(\Si, \O_\Si)$, and hence we can form
the standard symplectic disk bundle $(E_\Si, \w_{\can}) \to (\Si, \O_\Si)$.  This is the bundle from Biran's decomposition theorem \eqref{decom}, compare to \cite[Section 2]{Bir01} and \cite[Section 1]{Ops11}.

Returning to the standard symplectic disk bundle $\pi\co (E, \w) \to (\Si, \s)$ modeled on $P$, 
let $\th\co  [0,1] \to \R$ be nonnegative, zero in a neighborhood of one, and such that $\th(r)\co  E \to \R$ as a function of the radial coordinate $r$ is smooth.  We have an associated linear, order preserving map
\begin{equation}\label{Theta no F}
\T\co  C^\infty(\Si, \s) \to C^\infty(E, \w) \quad \mbox{defined as} \quad
H \mapsto \th \cdot\,\pi^{*}H.
\end{equation}
By integrating over the fiber, one can check that 
\begin{equation}\label{fiber}
	-\int_0^1 \th'(r)(1-r^2)^n\,dr \int_\Si H\, \s^{n-1} = \int_E \T(H)\, \w^n,
\end{equation}
so $\T$ preserves the property of a function having zero mean and hence induces a map
$\T\co  \P\H(\Si) \to \P\H(E)$.  At times we will be interested in the maps
$\T_\e$ which are induced by $\th_\e(r): [0,1] \to \R$, which smoothly interpolates between
$1-r^2$ when $r \leq 1-\e$ and zero when $r \geq 1-\tfrac{\e}{2}$.

For ease of notation and future applicability, will prove our lemmas in this general setting.
Note that the $\T$ from \eqref{Theta} in 
Theorems~\ref{thm: descent qm} and \ref{thm: descent} is just $\T$ from
\eqref{Theta no F} combined with the use of the symplectomorphism $F_\P$ from \eqref{decom} to view
$(E_\Si, \w_{\can}) = (E, \w)$ as an open subset of $M$.  So our lemmas will be applicable to
the setting of the theorems.

\subsection{An overview of the proofs of Theorem~\ref{thm: descent qm} and \ref{thm: descent}}
Before diving into the details of proving Theorem~\ref{thm: descent qm} and \ref{thm: descent},
let's step back and identify the main issue.  
We have a homogeneous quasi-morphism $\mu\co  \CH(E) \to \R$ and a map 
$\T\co  \PH(\Si) \to \PH(E)$, and we want to form a quasi-morphism
$\hat\mu\co  \CH(\Si) \to \R$.  Our situation can be summarized by the diagram: 
	\[
	\xymatrix{
	\PH(\Si) \ar[d]_{p_\Si} \ar[r]^{\T} & \PH(E) \ar[d] \ar[dr]^{\mu} \\
	\CH(\Si) \ar@{-->}@/_1.3pc/[rr]_{\hat\mu} & \CH(E) \ar[r]^{\quad\mu} & \R}
	\]
If $\T\co  \H(\Si) \to \H(E)$ preserved the Poisson brackets, then we would be done for $\T\co \PH(\Si) \to \PH(E)$ would be a homomorphism and would descend to a homomorphism $\widehat{\T}\co  \CH(\Si) \to \CH(E)$, which we could use to pullback $\mu$ to a quasi-morphism on $\CH(\Si)$.

Unfortunately $\T$ does not preserve the Poisson brackets, but since we only want 
a quasi-morphism, the full strength of a homomorphism is not necessary.  In Section~\ref{group}
we will formulate when a homogeneous quasi-morphism $\mu\co  G \to \R$ can be pulled back by a map $H \to G$ (Lemma~\ref{lem:pullback}), and when $\mu$ will
descend along a quotient $G \to H$ (Lemma~\ref{lem:pushforward}).  We will then show that 
$\mu\co  \PH(E) \to \R$ pulls back
along $\T_\e$ 
to a homogeneous quasi-morphism 
$\T_\e^*\mu\co  \PH(\Si) \to \R$, which then descends along $p_\Si$ to a homogeneous
quasi-morphism $\hat\mu\co  \CH(\Si) \to \R$. The technical heart of these two steps occupies Lemmas~\ref{lem: sgrad} and \ref{lem: homotopy}, from which the other parts of Theorem~\ref{thm: descent qm} and Theorem~\ref{thm: descent} follow.

\subsection{Two group theory lemmas}\label{group}

For these two lemmas, $G$ and $H$ will be groups and $\mu\co  G \to \R$ will be a homogeneous quasi-morphism with defect $D(\mu)$.

\begin{defn}
	A map $\vp\co  H \to G$ is \emph{$\mu$-quasi-homomorphism} if there is some $D(\vp, \mu)$
	such that
	\[
		\abs{\mu\param{\vp(h_1h_2)^{-1}\vp(h_1)\vp(h_2)}} \leq D(\vp, \mu)
		\quad\mbox{for all $h_1, h_2 \in H$}.
	\]
	$\vp$ is \emph{$\mu$-homogeneous} if there is some $C$ such that
	\[
		\abs{\mu\param{\vp(h^n)^{-1}\vp(h)^n}} \leq C \quad \mbox{for all $n \in \Z$ and $h \in H$}.
	\]
\end{defn}

\begin{lem}\label{lem:pullback}
	If $\vp\co H \to G$ is a $\mu$-quasi-homomorphism, then $\vp$ pulls $\mu$
	back to $\vp^*\mu \co  H \to \R$, a quasi-morphism on $H$,
	with defect bounded above by $D(\vp,\mu) + 2D(\mu)$.  If $\vp$ is also $\mu$-homogeneous,
	then $\vp^*\mu$ is a homogeneous quasi-morphism.
\end{lem}
\begin{proof}  We have that
	\begin{align*}
	\abs{\mu(\vp(h_1h_2)) \!-\! \mu(\vp(h_1)) \!-\! \mu(\vp(h_2))} &\leq
	\abs{-\mu(\vp(h_1h_2)) \!+\! \mu(\vp(h_1)\vp(h_2))} + D(\mu)\\
	&\leq \abs{\mu\param{\vp(h_1h_2)^{-1}\vp(h_1)\vp(h_2)}} + 2 D(\mu)\\
	&\leq D(\vp, \mu) + 2D(\mu),
	\end{align*}
	so $\vp^*\mu$ is a quasi-morphism.  If $\vp$ is $\mu$-homogeneous, then
	\begin{align*}
		\abs{\mu(\vp(h^n)) - \mu(\vp(h)^n)} &= \abs{\mu(\vp(h^n)^{-1}) + \mu(\vp(h)^n)}\\
		&\leq \abs{\mu\param{\vp(h^n)^{-1}\vp(h)^n}} + D(\mu) \leq C + D(\mu).
	\end{align*}
	Dividing through by $n$ and taking the limit, we see that $\vp^*\mu$ is its own homogenization.
\end{proof}

\begin{lem}\label{lem:pushforward}
	Let $\vp\co  G \to H$ be a surjective
	homomorphism such that $\mu$ is bounded on the kernel of $\vp$.  
	Then $\vp$ pushes $\mu$ forward to a homogeneous quasi-morphism
	\[
		\vp_*\mu\co H \to \R \quad \mbox{where} \quad (\vp_*\mu)(h) = \mu(g)
		\quad \mbox{for any $g \in \vp^{-1}(h)$},
	\]
	with defect $D(\mu)$.
\end{lem}
\begin{proof} The statement is clear provided that $\vp_*\mu$ is well defined.
Let $B \geq 0$ be such that $\abs{\mu(\ker \vp)} \leq B$.
Now if $g_1$ and $g_2$ both map to $h \in H$, then
$g_1^n$ and $g_2^n$ differ by an element of the kernel, so
\[
	n \abs{\mu(g_1) - \mu(g_2)} = \abs{\mu(g_1^n) - \mu(g_2^n)} \leq B + D(\mu).
\]
Now take the limit as $n$ goes to infinity, to get $\mu(g_1) = \mu(g_2)$.
\end{proof}

\subsection{The proofs of Theorem \ref{thm: descent qm} and Theorem~\ref{thm: descent}}\label{proofs}

We will need the following four lemmas, which are proved in Section \ref{lemma proofs}.
\begin{lem}\label{lem: sgrad}
	The map $\T\co  C^\infty(\Si, \s) \to C^\infty(E, \w)$ preserves the vanishing of Poisson brackets,
	namely for $H,K \in C^\infty(\Si)$ we have that
	\begin{equation}\label{poisson}
	\{\T(H),\T(K)\}_{\w}= \th(r)^2\{\pi^*H, \pi^*K\}_\w = \frac{\th(r)}{1-r^2}\, \T(\{H,K\}_{\s}).
	\end{equation}
	If $\{f_t\} \in \PH(\Si)$ and $\{\wt{f}_t\} \in \PH(E)$ are generated by
	$F \in \P\H(\Si)$ and $\T_\e(F) \in \P\H(E)$, then
	\begin{equation}\label{commute}
		\mbox{$f_t \circ \pi = \pi \circ \wt{f}_t$ when $r \leq 1-\e$.}
	\end{equation}
	The term measuring the failure of 
	$\T_\e\co  \P\H(\Si) \to \P\H(E)$ to be a homomorphism
	\begin{equation}\label{failure}
		\overline{\T_\e(F\#G)}\,\#\,(\T_\e(F)\#\T_\e(G)) \co  E \times [0,1] \to \R	
	\end{equation}
	vanishes at points in $E$ with $r \leq 1-\e$. This also holds for larger products as well,
	in particular for $\overline{\T_\e(F^{\#k})}\,\#\,(\T_\e(F)^{\#k})$.
\end{lem}
\begin{lem}\label{lem: displace}
	If $X\subset \Si$ is displaceable, then
	$\wt{X}_a = \pi^{-1}(X) \cap \{r \leq a\} \subset E$
	is displaceable if $a < 1$.
\end{lem}
\begin{lem}\label{lem: homotopy}
	Suppose that $F \in \P\H(\Si)$ generates a loop $f$, which is a null homotopic $[f] = \ind$ in 
	$\CH(\Si)$, and let $\T_\e(F) \in \P\H(E)$ generate the path $\wt{f}$.  Then as an element
	of $\CH(E)$, $[\wt{f}] = [\eta]$ where $\eta$ is generated by a normalized Hamiltonian
	on $E$ that vanishes when $r \leq 1-\e$.
\end{lem}
\begin{lem}\label{lem: scale}
Let $\mu\co  \CH(M) \to \R$ be a stable homogeneous quasi-morphism and let $H \in C^\infty(M)$
be a normalized Hamiltonian.  For any smooth function $\l\co  [0,1] \to \R$, 
\[
	\mu(\p_{\l H}) = \param{\int_0^1\l(t) dt}\mu(\p_H)
\]
where $\l H \in \P\H(M)$ is a time dependent normalized Hamiltonian.
\end{lem}

\subsubsection{The proof of Theorem~\ref{thm: descent qm}}
\begin{proof}[Proof that $\mu$ descends to $\hat\mu$.]
		By picking $\e_0$ small, we can take the neighborhood $U$ of $\Del_\P$
		to be the complement of $F_\P(r \leq 1-\e_0) \subset M$, and we will
		work with $\e < \e_0$.
		
		If $F, G \in \P\H(\Si)$ are normalized Hamiltonians, then by 
		\eqref{failure} in Lemma~\ref{lem: sgrad}
		\begin{equation}\label{ham1}
			\overline{\T_\e(F \# G)}\,\#\,(\T_\e(F) \# \T_\e(G)) \quad
			\mbox{and} \quad
			\overline{\T_\e(F^{\#k})}\,\#\, \T_\e(F)^{\#k}
		\end{equation}
		are normalized and supported on $U$.  Since $\mu$ restricted to
		$\CH_U(M)$ is the Calabi homomorphism, $\mu$ vanishes on \eqref{ham1}
		since they are normalized.
		Therefore by Lemma~\ref{lem:pullback}, $\T_\e$ pulls back $\mu$ to a 
		homogeneous quasi-morphism 	$\T_\e^*\mu$ on $\P\H(\Si)$ with defect at most 
		$2D(\mu)$.  It is independent of $\e$ since if $\e', \e < \e_0$, then
		$\mu$ vanishes on
		$\overline{\T_\e(F)}\, \#\, \T_{\e'}(F)$, for it is normalized and supported on $U$.
		Hence $\T_\e^*\mu = \T_{\e'}^*\mu$, for they are homogeneous quasi-morphism
		a bounded distance apart.
		
		If $F \in \P\H(\Si)$ generates a null-homotopic element
		in $\CH(\Si)$, then by Lemma~\ref{lem: homotopy}
		$\T_\e(F)$ generates an element
		$\eta \in \CH(M)$, which can also be generated by a normalized
		Hamiltonian in $\P\H(M)$ that is supported in $U$.
		Therefore by the Calabi property again, $\mu(\eta) = 0$, and hence
		the quasi-morphism $\T_\e^*\mu$ vanishes on the kernel of $\P\H(\Si) \to \CH(\Si)$.
		Therefore by Lemma~\ref{lem:pushforward}, $\T_\e^*\mu$ descends to
		$\CH(\Si)$ as the homogeneous quasi-morphism $\hat\mu$ in \eqref{qmorph descent}
		with defect at most $2D(\mu)$.	
\end{proof}

\begin{proof}[Proof that $\bar\mu$ inheriets the stability and Calabi properties from $\mu$] 
		Let $\mu$ be stable and have the Calabi property, and recall that
		$\bar\mu := \z(\th_\e)^{-1}\hat\mu$.
		Observe that for normalized functions $H, K\co  \Si \to \R$,
		\[
		\min_{\Si} (H - K) = \min_{M} (\T_\e(H) - \T_\e(K))
		\]
		and likewise for $\max$.  Hence the stability of $\bar\mu$ will follow
		from the stability of $\mu$.  The normalization constant $\z(\th_\e)^{-1}$ is independent
		of $\e$, for small $\e$, because $\th_\e$ and $\th_{\e'}$ Poisson commute and
		hence $\z(\th_\e) - \z(\th_{\e'}) = \z(\th_\e - \th_{\e'}) = 0$, for the difference is
		supported on the set $U$ where $\mu$ restricts to the Calabi homomorphism.
		
		To prove the Calabi property for $\bar\mu$, let $F\co \Si \times [0,1] \to \R$
		be supported on a displaceable set $V \subset \Si$ and let's denote
		$\vo = \vo(M) = \vo(\Si)$.  Consider the normalization terms for $F$ and $\T_\e(F)$:
		\[
			\l(t) = \vo^{-1} \int_\Si F_t\, \O_{\Si}^{n-1} \quad\mbox{and}\quad
			\eta(t) = \vo^{-1} \int_M \T_\e(F_t)\, \O^n = C_\e \l(t),
		\]
		where the last equality uses \eqref{fiber} and the notation 
		$C_\e = \vo^{-1}\int_M \th_\e\, \O^n$.  
		Since
		$\T_\e(\l) - \eta = \l\cdot(\th_\e - C_\e)$, by Lemma~\ref{lem: scale} we have that
		\[
		\mu(\p_{\T_\e(\l) - \eta}) = \param{\int_0^1 \l(t)\, dt} \mu(\p_{\th_\e - C_\e})
		= \Cal_V(F)\, \frac{\mu(\p_{\th_\e - C_\e})}{\vo}.
		\]
		By Lemma~\ref{lem: displace} the support of $\T_\e(F)$ is displaceable, so
		therefore 
		\[
		\mu(\p_{\T_\e(F) - \eta}) 
		= \int_0^1\int_M \T_\e(F_t)\, \O^n dt = C_\e \Cal_V(F)
		= \frac{\int_M \th_\e\, \O^n}{\vo}\, \Cal_V(F)
		\]
		Bringing this all together and using that $\T_\e(F)$ and $\T_\e(\l)$ commute gives	
		\begin{align*}
			\bar\mu(\p_{F-\l}) 
			&= \z(\th_\e)^{-1} 
			\param{\mu(\p_{\T_\e(F-\l)}) + \mu(\p_{\T_\e(\l) - \eta}) - \mu(\p_{\T_\e(\l) - \eta})}\\
			&= \z(\th_\e)^{-1}\param{\mu(\p_{\T_\e(F) - \eta}) - \mu(\p_{\T_\e(\l) - \eta})}\\
			&= \Cal_V(F)\, \z(\th_\e)^{-1} \frac{\int_M \th_\e\, \O^n - \mu(\p_{\th_\e - C_\e})}{\vo}
			= \Cal_V(F),
		\end{align*}
		as desired.		
	\end{proof}

\subsubsection{The proof of Theorem~\ref{thm: descent}}

\begin{proof}[Proof that $\z$ descends to $\bar\z_{\th}$]
	The monotonicity and normalization of $\bar{\z}_\th$ are immediate.
	As for quasi-linearity, if $F, G \in C^\infty(\Si)$ Poisson commute,
	then by \eqref{poisson}
	we have that 
	$\T(F), \T(G) \in C^\infty(M)$ Poisson commute as well.
	Therefore 
	\[
		\bar{\z}_\th(F + G) = \frac{\z(\T(F)+\T(G))}{\z(\th)} =
		\frac{\z(\T(F)) + \z(\T(G))}{\z(\th)} = \bar{\z}_\th(F) + \bar{\z}_\th(G),
	\]
	so we have that $\bar{\z}_\r$ is infact a symplectic quasi-state.
\end{proof}

\begin{proof}[Proof that $\bar\z_\th$ inherits properties from $\z$]
	
	Suppose that $\z$ vanishes on functions with displaceable support.
	If $H \in C^\infty(\Si)$ has displaceable support $X \subset \Si$,
	then by Lemma~\ref{lem: displace}, $\T(H)\co  M \to \R$ does as well
	and hence $\bar\z_{\th}(H) = \z(\T(H)) = 0$.
	
	Suppose that $\z\co  C^\infty(M) \to \R$ is $Ham(M)$ invariant.
	Let $\{f_t\}$ be a Hamiltonian isotopy generated by $F \in \P\H(\Si)$.  
	Let $\e$ be small enough so that $\th(r)$ vanishes when $r > 1- \e$, and let $\T_\e(F) \in \PH(M)$
	generate the Hamiltonian isotopy $\{\wt{f}_t\}$.
	On the image of $F_\P$ in $M$, by \eqref{commute} we have that
	$f_t \circ \pi \circ F_\P^{-1} = \pi \circ F_\P^{-1} \circ \wt{f}_t$, and hence
	for $H \in C^\infty(\Si)$ we have that
	\[
		\T(H \circ f_t)= \th \cdot (H \circ f_t \circ \pi^* \circ F_\P^{-1})
		=  (\th \cdot (H \circ \pi^* \circ F_\P^{-1}))\circ \wt{f}_t
		= \T(H) \circ \wt{f}_t.
	\]
	Therefore $\bar\z_\th$ is $Ham(\Si)$ invariant, for
	\[
	\bar\z_\th(H \circ f_t) = \frac{\z(\T(H \circ f_t))}{\z(\T(1))} =
	\frac{\z(\T(H) \circ \wt{f}_t)}{\z(\T(1))} = \bar\z_\th(H).
	\]	
	Finally, for the inequality \eqref{master} suppose that $\z$ satisfies it with
	$C(\z)$.  Then by \eqref{poisson},
	for $H, K \in C^\infty(\Si)$ we have that
	\begin{align*}
		\Pi_{\bar{\z}_\th}(H,K) &= \z(\th)^{-1}\, \Pi_\z(\T(H),\T(K))
		\leq \frac{C(\z)}{\z(\th)} \sqrt{\norm{\{\T(H),\T(K)\}}}\\
		&\leq \frac{C(\z)}{\z(\th)} \norm{\frac{\th(r)}{\sqrt{1-r^2}}} \sqrt{\norm{\{H,K\}}}.
	\end{align*}
	So therefore $C(\bar\z_\th) \leq \frac{C(\z)}{\z(\th)} \norm{\tfrac{\th(r)}{\sqrt{1-r^2}}}$.
\end{proof}

\subsection{The proofs of Lemmas~\ref{lem: sgrad}, \ref{lem: displace}, \ref{lem: homotopy},
and \ref{lem: scale}}\label{lemma proofs}

\begin{proof}[Proof of Lemma~\ref{lem: sgrad}]
	We will start by computing $\pi_*(\sgrad \pi^*H)$.  For a vector $X \in T_{(p,r)}E$
	\begin{align*}
		-dH(\pi_*X) &= \w(\sgrad \pi^*H, X)\\
		&= (\pi^*\s + r^2 d\a)(\sgrad \pi^*H, X) &\mbox{by
		$\i_{\sgrad \pi^*H}(dr \wedge \a) = 0$}\\
		&= (1 - r^2)\, \s(\pi_*(\sgrad \pi^*H), \pi_*X) &\mbox{by $\s = -\pi_*d\a$}
	\end{align*}
	so we conclude that $\pi_*(\sgrad\pi^*H) = \frac{1}{1-r^2}\,\sgrad H$.  Therefore
	\begin{equation}\label{push sgrad theta}
		\pi_*(\sgrad\T(H)) = \th(r)\,\pi_*(\sgrad\pi^*H) = \frac{\th(r)}{1-r^2}\,\sgrad H.
	\end{equation}
	Now since $\pi^*H$ and $\th(r)$ Poisson commute, the computation for
	\eqref{poisson} is the following:
	\begin{align*}
		\{\T(H), \T(K)\}_\w &=
		\th(r)^2\{\pi^*H, \pi^*K\}_\w = \th(r)^2\, dH(\pi_*\sgrad \pi^*K)\\
		&= \frac{\th(r)^2}{1-r^2}\, \pi^*(\{H,K\}_\s) = \frac{\th(r)}{1-r^2}\, \T(\{H,K\}_{\s}).
	\end{align*}
	
	Going to the path space side, using \eqref{push sgrad theta} 
	we see that $\pi\circ \wt{f}_t\co  E \to \Si$ satisfies the 
	differential equation
	\[
		\d_t(\pi \circ \wt{f}_t)
		= \param{\frac{\th_\e(r)}{1-r^2} \sgrad F_t}_{\pi \circ \wt{f}_t}
		\quad \mbox{with initial condition } \pi \circ \wt{f}_0 = \pi.
	\]
	Staring at a point with $\th_\e(r) = 1-r^2$, then $f_t \circ \pi$ also satisfies
	this equation, so \eqref{commute} follows by the uniqueness of solutions to ODEs.
	Denoting by $\wt{(fg)}_t$ the Hamiltonian path generated by 
	$\T_\e(F \# G)$ we have that
	\begin{align*}
		\overline{\T_\e(F\#G)}\,\#\,(\T_\e(F)\#\T_\e(G)) 
		&= (- \T_\e(F\#G) + \T_\e(F)\#\T_\e(G))\circ \wt{(fg)}_t\\
		&= \th_\e(r) \,\param{-G \circ f_t^{-1} \circ \pi + G \circ \pi \circ \wt{f}_t^{-1}} \circ \wt{(fg)}_t,
	\end{align*}
	which vanishes when $r \leq 1-\e$ since then we can use that 
	$f_t^{-1} \circ \pi = \pi \circ \wt{f}^{-1}_t$.  The proof for larger products is the same.	
\end{proof}

\begin{proof}[Proof of Lemma~\ref{lem: displace}]
	Let $F \in \P\H(\Si)$ generate the isotopy $f$ which displaces $X$.  Pick $\e > 0$ so that
	$1-\e > a$ and let $\wt{f}$ be the isotopy generated by $\T_\e(F) \in \P\H(E)$.  It follows
	from \eqref{commute} that $\wt{f}$ displaces $\wt{X}_a$. 
\end{proof}

\begin{proof}[Proof of Lemma~\ref{lem: homotopy}]
	Let $\vp_t^s$ be a homotopy of loops in $Ham(\Si)$ based at $\ind$, between the
	loop $\{\vp^0_t = f_t\}_{t\in [0,1]}$ and the constant loop $\{\vp^1_t = \ind\}_{t\in [0,1]}$.  
	For $s$ fixed, let $F_t^s\co  \Si \to \R$ be the 
	Hamiltonian in $\P\H(\Si)$ generating
	the Hamiltonian loop $\{\vp^s_t\}_{t\in [0,1]}$ in $Ham(\Si)$, via
	\[
		\d_t \vp_t^s = (\sgrad F_t^s)_{\vp_t^s} \quad
		\mbox{with}\quad \vp_0^s = \ind.
	\]
	Note that $F_t^0 = F_t$ and $F_t^1 = 0$.
	While for $t$ fixed, let $G_t^s\co  \Si \to \R$ be the
	Hamiltonian in $\P\H(\Si)$ generating
	the Hamiltonian path $\{\vp_t^s\}_{s \in [0,1]}$ in $Ham(\Si)$, via
	\[
		\d_s \vp_t^s = (\sgrad G_t^s)_{\vp_t^s} \quad\mbox{with}\quad
		\vp_t^0 = f_t.
	\]
They are related by the equation \cite[Proposition I.1.1]{Ban78}
\begin{equation}\label{dhf1}
	\d_s F_t^s = \d_t G_t^s + \{F_t^s, G_t^s\}.
\end{equation}

Fixing $s$,
the Hamiltonian $\T_\e(F^s_t)\co E \to \R$ in $\P\H(E)$ will generate a Hamiltonian path
$\{\psi^s_t\}_{t\in [0,1]}$ in $Ham(E)$, via
\[
	\d_t \psi_t^s = (\sgrad \T_\e(F_t^s))_{\psi_t^s} \quad\mbox{with}
	\quad \psi_0^s = \ind.
\]
As $s$ varies, $\psi_s^t$ will be a homotopy of Hamiltonian paths in $Ham(E)$, between 
the paths $\{\psi_t^0 = \wt{f}_t\}_{t\in [0,1]}$ and $\{\psi_t^1 = \ind\}_{t\in [0,1]}$.
However this will not be a homotopy of loops, since in particular 
$\{\wt{f}_t\}_{t\in [0,1]}$ may not be a loop.
Letting $t$ be fixed,
let $H_t^s\co  E \to \R$ be the
Hamiltonian in $\P\H(E)$ generating the Hamiltonian path $\{\psi_t^s\}_{s \in [0,1]}$ in $Ham(E)$, via
\[
	\d_s \psi_t^s = (\sgrad H_t^s)_{\psi_t^s} \quad \mbox{with}
	\quad \psi_t^0 = \wt{f}_t.
\]
Just as in \eqref{dhf1}, we have the relation
\begin{equation}\label{dhf2}
	\d_s \T_\e(F_t^s) = \d_t H_t^s + \{\T_\e(F_t^s), H_t^s\}.
\end{equation}

Define the path $\{\eta_u = (\psi_1^u)^{-1}\psi_1^0\}_{u \in [0,1]}$ in $Ham(E)$, which
is a path from $\eta_0=\ind$ to $\eta_1 = \psi_1^0$, using that $\psi_t^1 = \ind$. 
Observe that $\Psi_t^s = \psi_t^s\eta_{st}$ is a homotopy
of paths in $Ham(E)$, between the paths
\[
	\{\Psi^0_t = \psi^0_t = \tilde{f}_t\}_{t \in [0,1]}
	\quad \mbox{and} \quad \{\Psi_t^1 = \eta_t\}_{t \in [0,1]}.
\]
Since $\Psi$ is a homotopy of paths with fixed endpoints,
\[\Psi^s_0 = \ind \quad \mbox{and} \quad \Psi_1^s = \psi_1^s \eta_s = \psi_1^0,
\]
we have proved that $[\wt{f}] = [\eta]$ in $\CH(E)$. 

The path $\eta$ is generated by the normalized Hamiltonian
$(-H_1^u \circ \psi_1^u)$, which we claim vanish at points with $r \leq 1-\e$.
It suffices to prove this for $H_1^u$ since $\psi_t^s$ preserves the level sets of $r$.
Applying $\T_\e$ to \eqref{dhf1} and using \eqref{poisson} of Lemma~\ref{lem: sgrad} gives that
	\begin{equation}\label{dhf3}
		\d_s \T_\e(F_t^s) = \d_t \T_\e(G_t^s) + \frac{(1-r^2)}{\th_\e(r)}\, \{\T_\e(F_t^s), \T_\e(G_t^s)\}.
	\end{equation}
	When $r \leq 1-\e$, we have that $\th_\e(r) = 1-r^2$ and so the differential
	equations \eqref{dhf2} and \eqref{dhf3} agree.  Therefore by the method of characteristics for PDE's,
	$H_t^s = \T_\e(G_t^s)$ when $r \leq 1 - \e$.
	The fact that $\vp$ is a homotopy of paths with fixed endpoints means that
	$G_1^s = 0$, and therefore $H_1^s = 0$ when $r \leq 1-\e$.
\end{proof}

\begin{proof}[Proof of Lemma~\ref{lem: scale}]
If $L(t) = \int_0^t \l(s)ds$ and $h_t$ is the flow generated by $\sgrad(H)$, then
the Hamiltonian $\l H$ generates the Hamiltonian isotopy $\p_{\l H} = \{h_{L(t)}\}_{t\in [0,1]}$.
By a time reparameterization, in $\CH(M)$ we have that $[\{h_{L(t)}\}] = [\{h_{tL(1)}\}]$ and
hence $\mu(\p_{\l H}) = \mu(\p_{L(1)H})$.  If $m$ is an integer, then
$\mu(\p_{mH}) = \mu(\p_H^m) = m\mu(\p_H)$ since $\mu$ is homogeneous and
hence $\mu(\p_{a H}) = a \mu(\p_{H})$ when $a$ is rational.  Since $\mu$ is stable,
this extends to any real scalar, so in particular $\mu(\p_{L(1)H}) = L(1) \mu(\p_{H})$.
\end{proof}

\subsection{The proof of Theorem~\ref{thm: stably}}\label{proofs2}
\begin{proof}[Proof of Theorem~\ref{thm: stably}]
	Since $X \times S^1$ is displaceable in $M \times T^*S^1$,
	it is also displaceable in $M \times S^1 \times [-R, R]$ where the last coordinate is
	the vertical direction in $T^*S^1$.  Capping $S^1 \times [-R, R]$ off on the top and bottom
	with half spheres creates an sphere.  
	Therefore by Moser we may assume that $X \times S^1 \subset M \times (S^2, \O)$ is
	displaceable, where $(S^2, \Omega)$ is a round sphere in $\R^3$ 
	with coordinates $(x,y,z)$, $\O$ is the induced area form, 
	and $S^1$ is the equator $\{z=0\} \subset S^2$.

	Let $U$ and $V = \{\abs{z} < 2\e\}$ be open neighborhoods of $X$ and $S^1$
	respectively such that $U \times V \subset M \times S^2$ is still displaceable.  
	Let $\rho\co  M \to \R$ be a cut-off function 
	with support in $U$ that is a constant $1$ near $X$.  On $S^2$, let $V_1 = \{z> \e\}$ and
	$V_2 = \{z< -\e\}$, so $\{V=V_0, V_1, V_2\}$ is an open cover of $S^2$
	and let $\phi_i(z)$ be a Poisson commuting subordinate partition of unity. 
	Note that each $U \times V_i \subset M \times S^2$ is displaceable.
	
	Let $\pi_1\co M\times S^2 \to M$ and $\pi_2\co  M \times S^2 \to S^2$ be the projections
	and let $\z(a\otimes[S^2], \cdot)$ be the partial symplectic quasi-state on 
	$M \times S^2$ associated to the
	idempotent $a \otimes [S^2] \in QH_{2n+2}(M \times S^2)$, see \cite[Section 3.5]{EntPol09RS}.  Since
	the $\phi_i$ pullback to a Poisson commuting partition on unity on $M \times S^2$,
	\[
		\z(a \otimes [S^2], \pi_1^* \rho) 
		= \z\left(a \otimes [S^2],\,\sum_{i=0}^2 \pi_1^*\rho \cdot \pi_2^* \phi_i\right)
		= 0.
	\]
	The above vanishes by the partial additivity and vanishing property
	\cite[Theorem 3.6]{EntPol09RS} of partial symplectic quasi-states,
	using that the $\pi_1^*\rho \cdot \pi_2^* \phi_i$ Poisson commute and have displaceable
	support.
	Finally by the product formula for spectral invariants \cite[Theorem 5.1]{EntPol09RS}, we have that
	\[
		\z(a, \rho) = \z(a \otimes [S^2], \pi_1^*\rho) = 0.
	\]
	It now follows that from the definition of $\tau(a)$ that $\tau(a, X) = 0$.
\end{proof}

\section{Questions}

We will end this paper with two questions regarding 
the process of passing from a quasi-morphism $\mu: \CH(M) \to \R$ to its reduction 
$\bar\mu: \CH(\Si) \to \R$ as described in Theorem~\ref{thm: descent qm}.
\begin{enumerate}
\item If $\mu$ descends from the universal cover to $\mu: Ham(M) \to \R$, does its reduction also descend to $\bar\mu: Ham(\Si) \to \R$? 
\item If $\mu = \mu(a, \cdot):  \CH(M) \to \R$ is a spectral quasi-morphism for some idempotent
$a \in QH_{2n}(M)$, then
does its reduction $\bar\mu = \mu(\bar{a}, \cdot)$ come for some
idempotent $\bar{a} \in QH_{2n-2}(\Si)$?
\end{enumerate}
The trouble with proving 1) is that $\T_{\e}: \PH(M) \to \PH(\Si)$ does not preserve loops due to the
discrepancy between $\th_{\e}(r)$ and $1-r^{2}$.  
For 2) the trouble is that there is no natural candidate for $\bar{a}$ since there is no known natural map 
$QH_{2n}(M) \to QH_{2n-2}(\Si)$ preserving the ring structure.


\bibliographystyle{alpha}
\bibliography{symplectic}

\begin{thebibliography}{FOOO11}

\bibitem[Aar91]{Aar91}
J.~F. Aarnes.
\newblock Quasi-states and quasi-measures.
\newblock {\em Adv. Math.}, 86(1):41--67, 1991.

\bibitem[Ban78]{Ban78}
A.~Banyaga.
\newblock Sur la structure du groupe des diff{\'e}omorphismes qui
  pr{\'e}servent une forme symplectique.
\newblock {\em Comment. Math. Helv.}, 53(2):174--227, 1978.

\bibitem[BC01]{BirCie01}
P.~Biran and K.~Cieliebak.
\newblock Symplectic topology on subcritical manifolds.
\newblock {\em Comment. Math. Helv.}, 76(4):712--753, 2001.

\bibitem[BEP11]{BuhEntPol11}
L.~Buhovsky, M.~Entov, and L.~Polterovich.
\newblock Poisson brackets and symplectic invariants.
\newblock arXiv:1103.3198v1, 2011.

\bibitem[Bir01]{Bir01}
P.~Biran.
\newblock Lagrangian barriers and symplectic embeddings.
\newblock {\em Geom. Funct. Anal.}, 11(3):407--464, 2001.

\bibitem[Bir06]{Bir06}
P.~Biran.
\newblock Lagrangian non-intersections.
\newblock {\em Geom. Funct. Anal.}, 16(2):279--326, 2006.

\bibitem[BJ11]{BirJer11}
P.~Biran and Y.~Jerby.
\newblock The symplectic topology of projective manifolds with small dual.
\newblock Preprint, 2011.

\bibitem[BK11]{BirKha11}
P.~Biran and M.~Khanevsky.
\newblock A {F}loer-{G}ysin exact sequence for {L}agrangian submanifolds.
\newblock arXiv:1101.0946v1, 2011.

\bibitem[Bor11]{Bor11a}
M.~S. Borman.
\newblock Quasi-states, quasi-morphisms, and the moment map.
\newblock arXiv:1105.1805, 2011.

\bibitem[Buh10]{Buh10}
L.~Buhovsky.
\newblock The {$2/3$}-convergence rate for the {P}oisson bracket.
\newblock {\em Geom. Funct. Anal.}, 19(6):1620--1649, 2010.

\bibitem[Cal09]{Cal09}
D.~Calegari.
\newblock {\em scl}, volume~20 of {\em MSJ Memoirs}.
\newblock Mathematical Society of Japan, Tokyo, 2009.

\bibitem[CV08]{CarVit08}
F.~Cardin and C.~Viterbo.
\newblock Commuting {H}amiltonians and {H}amilton-{J}acobi multi-time
  equations.
\newblock {\em Duke Math. J.}, 144(2):235--284, 2008.

\bibitem[Don96]{Don96}
S.~K. Donaldson.
\newblock Symplectic submanifolds and almost-complex geometry.
\newblock {\em J. Differential Geom.}, 44(4):666--705, 1996.

\bibitem[EP03]{EntPol03}
M.~Entov and L.~Polterovich.
\newblock Calabi quasimorphism and quantum homology.
\newblock {\em Int. Math. Res. Not.}, (30):1635--1676, 2003.

\bibitem[EP06]{EntPol06}
M.~Entov and L.~Polterovich.
\newblock Quasi-states and symplectic intersections.
\newblock {\em Comment. Math. Helv.}, 81(1):75--99, 2006.

\bibitem[EP08]{EntPol08}
M.~Entov and L.~Polterovich.
\newblock Symplectic quasi-states and semi-simplicity of quantum homology.
\newblock In {\em Toric topology}, volume 460 of {\em Contemp. Math.}, pages
  47--70. Amer. Math. Soc., Providence, RI, 2008.

\bibitem[EP09]{EntPol09RS}
M.~Entov and L.~Polterovich.
\newblock Rigid subsets of symplectic manifolds.
\newblock {\em Compos. Math.}, 145(3):773--826, 2009.

\bibitem[EP10a]{EliPol10}
Y.~Eliashberg and L.~Polterovich.
\newblock Symplectic quasi-states on the quadric surface and {L}agrangian
  submanifolds.
\newblock arXiv:1006.2501v1, 2010.

\bibitem[EP10b]{EntPol10}
M.~Entov and L.~Polterovich.
\newblock {$C^0$}-rigidity of {P}oisson brackets.
\newblock In {\em Symplectic topology and measure preserving dynamical
  systems}, volume 512 of {\em Contemp. Math.}, pages 25--32. Amer. Math. Soc.,
  Providence, RI, 2010.

\bibitem[EPZ07]{EntPolZap07}
M.~Entov, L.~Polterovich, and F.~Zapolsky.
\newblock Quasi-morphisms and the {P}oisson bracket.
\newblock {\em Pure Appl. Math. Q.}, 3(4, Special Issue: In honor of Grigory
  Margulis. Part 1):1037--1055, 2007.

\bibitem[FOOO11]{FukOhOht11a}
K.~Fukaya, Y.-G. Oh, H.~Ohta, and K.~Ono.
\newblock Spectral invariants with bulk quasimorphisms and {L}agrangian {F}loer
  theory.
\newblock arXiv:1105.5123v1, 2011.

\bibitem[Gir02]{Gir02}
E.~Giroux.
\newblock G{\'e}om{\'e}trie de contact: de la dimension trois vers les
  dimensions sup{\'e}rieures.
\newblock In {\em Proceedings of the {I}nternational {C}ongress of
  {M}athematicians, {V}ol. {II} ({B}eijing, 2002)}, pages 405--414, Beijing,
  2002. Higher Ed. Press.

\bibitem[G{\"u}r08]{Gur08}
B.~G{\"u}rel.
\newblock Totally non-coisotropic displacement and its applications to
  {H}amiltonian dynamics.
\newblock {\em Commun. Contemp. Math.}, 10(6):1103--1128, 2008.

\bibitem[Kot04]{Kot04}
D.~Kotschick.
\newblock What is{$\dots$}a quasi-morphism?
\newblock {\em Notices Amer. Math. Soc.}, 51(2):208--209, 2004.

\bibitem[LM96]{LalMcD96}
F.~Lalonde and D.~McDuff.
\newblock The classification of ruled symplectic {$4$}-manifolds.
\newblock {\em Math. Res. Lett.}, 3(6):769--778, 1996.

\bibitem[MS98]{McDSal98}
D.~McDuff and D.~Salamon.
\newblock {\em Introduction to symplectic topology}.
\newblock Oxford Mathematical Monographs. The Clarendon Press Oxford University
  Press, New York, second edition, 1998.

\bibitem[Ops11]{Ops11}
E.~Opshtein.
\newblock Polarizations and symplectic isotopies.
\newblock {\em J. Symplectic Geom. (to appear)}, 2011.
\newblock arXiv:0911.3601v1.

\bibitem[Ost06]{Ost06}
Y.~Ostrover.
\newblock Calabi quasi-morphisms for some non-monotone symplectic manifolds.
\newblock {\em Algebr. Geom. Topol.}, 6:405--434 (electronic), 2006.

\bibitem[Pol01]{Pol01}
L.~Polterovich.
\newblock {\em The geometry of the group of symplectic diffeomorphisms}.
\newblock Lectures in Mathematics ETH Z{\"u}rich. Birkh\"auser Verlag, Basel,
  2001.

\bibitem[She11]{She11}
E.~Shelukhin.
\newblock The action homomorphism, quasimorphisms and moment maps on the space
  of compatible almost complex structures.
\newblock arXiv:1105.5814v1, 2011.

\bibitem[Tev03]{Tev03}
E.~A. Tevelev.
\newblock Projectively dual varieties.
\newblock {\em J. Math. Sci. (N. Y.)}, 117(6):4585--4732, 2003.

\bibitem[Ush10]{Ush10}
M.~Usher.
\newblock Deformed {H}amiltonian {F}loer theory, capacity estimates, and
  {C}alabi quasimorphisms.
\newblock arXiv:1006.5390v1, 2010.

\end{thebibliography}


\end{document}